\documentclass[12pt]{article}
\usepackage{babel}
\usepackage{isolatin1,times,a4wide,doublespace,amsmath,latexsym,amssymb}
\newtheorem{theorem}{Theorem}[section]
\newtheorem{lemma}[theorem]{Lemma}
\newtheorem{proposition}[theorem]{Proposition}

\newcommand{\tr}{\mathrm{tr}}

\newcommand{\Span}{\mathrm{span}}
\newcommand{\codim}{\mathrm{codim}}

\newcommand{\GL}{\mathrm{GL}}
\newcommand{\GAP}{\mathrm{GAP}}

\newtheorem{definition}[theorem]{Definition}

\newtheorem{corollary}[theorem]{Corollary}

\newcommand{\bproof}{\noindent{\bf Proof: }}
\newcommand{\eproof}{\hfill $\Box$\\}
\newcommand{\bremark}{\noindent{\bf Remark: }}
\newcommand{\eremark}{\hfill \\}
\newcommand{\w}{\omega}

\newcommand{\e}{\varepsilon}

\newcommand{\cN}{{\cal N}}

\newcommand{\cL}{{\cal L}}
\newcommand{\cS}{{\cal S}}
\newcommand{\cB}{{\cal B}}

\newcommand{\bN}{{\mathbb N}}
\newcommand{\bR}{{\mathbb R}}

\begin{document}
\title{The solution to the Maurey extension problem for Banach spaces
with the Gordon-Lewis property and related structures}
 
\author{P.G.\ Casazza\footnote{Supported by NSF grant DMS 970618.} \and
N.J.\ Nielsen\footnote{supported by the Danish Natural Science Research
Council, grants 9503296 and 9801867.}}
\date{ }
\maketitle
 
\begin{abstract}
The main result of this paper states that if a Banach space $X$ has the property that every
bounded operator from an arbitrary subspace of $X$ into an arbitrary
Banach space of cotype 2 extends to a bounded operator on $X$, then
$B(\ell_{\infty},X^*)=\Pi_2(\ell_{\infty},X^*)$. If in addition $X$ has
the Gaussian average property, then it is of type 2. This implies that
the same conclusion holds if $X$ has the Gordon-Lewis property (in
particular $X$ could be a Banach lattice) or if $X$ is isomorphic to a
subspace of a Banach lattice of finite cotype, thus solving the Maurey
extension property for these classes of spaces.

The paper also contains a detailed study of the property of extending
operators with values in $\ell_p$-spaces, $1\le p<\infty$.

\end{abstract}

\section*{Introduction}

In 1974 Maurey \cite{M1} proved that if $X$ is a Banach space of type 2,
then every bounded operator from an arbitrary subspace of $X$ to an
arbitrary Banach space $Y$ of cotype 2 admits a bounded extension from
$X$ to $Y$. Since then it has been an open problem whether this property
known as the Maurey extension property characterizes Banach spaces of
type 2. Since it follows from \cite{MIP} that a Banach space with this
property is of weak type 2, the answer to the problem is clearly
affirmative for the class of spaces where weak type 2 is equivalent to
type 2, e.g. rearrangement invariant function spaces.

The main result of this paper states that if a Banach space $X$ has the Maurey
extension property, then every bounded operator from an
$L_{\infty}$-space to $X^*$ is 2-summing. If in addition $X$ has Gaussian
average property $\GAP$ (as defined in \cite{CN1}), then it is of type
2. This implies that the answer to the problem is also affirmative for 
Banach spaces which have the Gordon-Lewis property, in particular Banach lattices, as well as for
Banach spaces which are isomorphic to subspaces of Banach lattices of
finite cotype.

It is not known in general whether the condition
$B(\ell_{\infty},X^*)=\Pi_2(\ell_{\infty},X^*)$ implies that $X^*$ is of
cotype 2 or equivalently in the case above that $X$ is of type 2. It
seems at the moment that $\GAP$ is the weakest known condition to ensure
this for K-convex spaces. It should be noted that every space of type 2
has $\GAP$.

We shall say that a Banach space $X$ has $M_p$, $1\le p<\infty$, if
every bounded operator from a subspace of $X$ to $\ell_p$ admits a bounded
extension to $X$. Another major result of the paper states that $M_p$,
$2<p<\infty$, characterizes Hilbert spaces among Köthe function spaces on
$[0,1]$. Finally we investigate $M_p$, $1\le p\le 2$ in detail and
prove that $M_1$ is equivalent to $M_p$, $1<p<2$ and that $M_1$ implies
$M_2$.

It is an open problem whether $M_2$ implies $M_1$ and whether $M_1$ or
$M_2$ imply the Maurey extension property.

We now wish to discuss the arrangement of this paper in greater detail.

In Section 1 of the paper we prove some general results on extensions of
operators which are needed to prove the main results. Some of them are
probably of interest in their own right. Section 2 is devoted to the main results stated
above while Section 3 contains the investigation of the properties
$M_p$, $1\le p\le 2$, and the proof of the implications
$M_1\Leftrightarrow M_p$, $1<p<2$, and $M_1\Rightarrow M_2$. 

\section*{Acknowledgement}

The authors are indebted to Nigel Kalton who drew our attention to the
spaces $\ell_p(\delta,2)$, \mbox{$2<p<\infty$} in order to prove that $M_p$ does
not have $M_r$ for $2<p<r<\infty$. This subsequently lead to the idea of
the proof of our main result.

Spaces like $\ell_p(\delta,2)$ were first considered by Rosenthal in his
construction of new $\cL_p$-spaces \cite{R}.
  
\setcounter{section}{-1} 
\section{Notation and Preliminaries}
\label{sec0}

In this paper we shall use the notation and terminology commonly used in
Banach space theory as it appears in \cite{LT1}, \cite{LT2} and
\cite{T}. $B_X$ shall always denote the closed unit ball of the Banach
space $X$.

If $X$ and $Y$ are Banach spaces, then $B(X,Y)$ ($B(X)=B(X,X)$) denotes the
space of all bounded linear operators from $X$ to $Y$ and throughout the
paper we shall identify $X\otimes Y$ with the space of all
$\w^*$-continuous finite rank operators from $X^*$ to $Y$ in the
canonical manner. Further if $1\le p<\infty$, we let $\pi_p(X,Y)$ denote
the space of all $p$-summing operators from $X$ to $Y$ equipped with the
$p$-summing norm $\pi_p$; $I_p(X,Y)$ denotes the space of all strictly
$p$-integral operators from $X$ to $Y$ equipped with the strict $p$-integral
norm $i_p$ and $N_p(X,Y)$ denotes the space of all $p$-nuclear operators
from $X$ to $Y$ equipped with the $p$-nuclear norm
$\nu_p$. $X\otimes_\pi Y$ denotes the completion of $X\otimes Y$ under
the largest tensor norm $\pi$ on $X\otimes Y$.

We recall that if $1\le p\le\infty$, then an operator $T\in B(X,Y)$ is
said to factor through $L_p$ if it admits a factorization $T=BA$ where
$A\in B(X,L_p(\mu))$ and $B\in B(L_p(\mu),Y)$ for some measure $\mu$
and we denote the space of all operators which factor through $L_p$ by
$\Gamma_p(X,Y)$. If $T\in\Gamma_p(X,Y)$, then we define
\[
\gamma_p(X,Y) = \inf\{\|A\| \|B\| \mid T=BA,\quad  \mbox{$A$ and $B$ as
above}\};
\]
$\gamma_p$ is a norm on $\Gamma_p(X,Y)$ turning it into a Banach
space. All these spaces are operator ideals and we refer to the above
mentioned books, \cite{DJT} and \cite{K1} for further details.

In the formulas of this paper we shall, as is customary, interpret
$\pi_\infty$ as the operator norm and $i_\infty$ as the
$\gamma_\infty$-norm.

We let $(r_n)$ denote the sequence of Rademacher functions on $[0,1]$
and recall that a Banach space $X$ is said to be of type p, $1\le p \le
2$ (respectively cotype p, $2\le p < \infty)$, if there is a constant $K\ge 1$ so that for all finite sets
$\{x_1,x_2,\dots,x_n\}\subseteq X$ we have
\begin{equation}
\label{eqa}   
\big(\int^1_0 \big\|\sum_{j=1}^n r_j(t)x_j\big\|^p
dt\big)^\frac1p \le K \big(\sum_{j=1}^n \|x_j\|^p \big)^\frac1p
\end{equation}
(respectively
\begin{equation}
\label{eqb}
\big(\sum_{j=1}^n \|x_j\|^p \big)^\frac1p \le\big(\int^1_0 \big\|\sum_{j=1}^n r_j(t)x_j\big\|^pdt\big)^\frac1p).
\end{equation}

The smallest constant $K$ which can be used in (\ref{eqa})
(respectively (\ref{eqb})) is denoted by $K^p(X)$ (respectively
$K_p(X)$). 

A Banach space $X$ is said to be of weak type 2 if there is a constant
$C$ and a $\delta$, $0<\delta<1$, so that whenever $E\subseteq X$ is a
subspace, $n\in\bN$ and $T\in B(E,\ell_2^n)$, then there is an orthogonal
projection $P$ on $\ell_2^n$ of rank larger than $\delta n$ and an
operator $S\in B(X,\ell_2^n)$ with $Sx=PTx$ for all $x\in E$ and
$\|S\|\le C\|T\|$.

Similarly $X$ is called a weak cotype 2 if there is a constant $C$ and a
$\delta$, $0<\delta<1$, so that whenever $E\subseteq X$ is a finite
dimensional subspace, then there is a subspace $F\subseteq E$ so that
$\dim F\ge \delta \dim E$ and $d(F,\ell_2^{\dim F})\le C$.

Our definitions of weak type 2 and weak cotype 2 space are not the
original ones, but are chosen out of the many equivalent
characterizations given by Pisier \cite{P2}.

Following \cite{GJN} we shall say that a Banach space $X$ has
$\GL(p,q)$, $1\le p,q\le\infty$, if there is a constant $K$ so that for
all Banach spaces $Y$ and all $T\in X^*\otimes Y$ we have $i_q(T)\le
K\pi_p(T^*)$. The smallest constant $K$ which can be used in this
inequality is denoted by $\GL_{p,q}(X)$. We note that $\GL(1,\infty)$
corresponds to the classical Gordon-Lewis property $\GL$ see
\cite{GL}. $X$ is said to have the Gordon -Lewis property $\GL_2$ if
every 1-summing operator from $X$ to a Hilbert space factors through an $L_1$-space.

If $n\in\bN$ and $T\in B(\ell_2^n,X)$, then following \cite[\S12]{T} we
define the $\ell$-norm of $T$ by
\[
\ell(T) = \big(\int_{\ell_2^n}\|Tx\|^2d\gamma(x)\big)^\frac12
\]
where $\gamma$ is the canonical Gaussian probability measure on $\ell_2^n$.

A Banach space $X$ is said to have the Gaussian Average Property
(abbreviated $\GAP$) \cite{CN1} if there is a constant $K$ so that
$\ell(T)\le K\pi_1(T^*)$ for every $T\in B(\ell^n_2,X)$ and every $n\in\bN$.

We shall also need some notation on subspaces of Banach lattices and on operators with ranges in a Banach
lattice. Recall that if $X$ is a Banach space and $L$ is a Banach
lattice, then an operator $T\in B(X,L)$ is called order bounded
\cite{NJN1} if there exists a $z\in L$, $z\ge 0$ so that
\begin{equation}
\label{eq0.1}
|Tx| \le \|x\| z\quad\mbox{for all $x\in X$}
\end{equation}
and the order bounded norm $\|T\|_m$ is defined by
\begin{equation}
\label{eq0.2}
\|T\|_m = \inf\{\|z\| \mid z\ \mbox{can be used in \eqref{eq0.1}}\}.
\end{equation}
We let $\cB(X,L)$ denote the space of all order bounded operators from
$X$ to $L$ equipped with the norm $\|\cdot\|_m$. It is readily seen to
be a Banach space and a left ideal. $X^*\otimes_m L$ shall denote the
closure of $X^*\otimes L$ in $\cB(X,L)$ under the norm $\|\cdot\|_m$.

If $X$ be a subspace of a Banach lattice $L$ and $1\le p< \infty$, then
we shall say that $X$ is $p$-convex in $L$ (respectively $p$-concave in $L$) if there
is a constant $K\ge 1$ so that for all finite sets
$\{x_1,x_2,\dots,x_n\}\subseteq X$ we have 
\[
\|(\sum^n_{j=1} |x_j|^p)^{\frac1p}\| \le K(\sum^n_{j=1}
\|x_j\|^p)^{\frac1p}
\]
(respectively
\[
(\sum^n_{j=1} \|x_j\|^p)^{\frac1p} \le K\|(\sum^n_{j=1}
|x_j|^p)^{\frac1p}\| ).
\]
Note that these inequalities depend on the embedding of $X$ into
$L$. $L$ is called $p$-convex (respectively $q$-concave) if the above
inequalties hold for every finite set of vectors in $L$.

If $E$ is a Banach space and $T\in B(E,X)$, then $T$ is called $p$-convex
if there exists a constant $K\ge 0$ so that for all finite sets
$\{x_1,x_2,\dots,x_n\}\subseteq E$ we have
\[
\|(\sum^n_{j=1} |Tx_j|^p)^{\frac1p}\| \le K(\sum^n_{j=1}
\|x_j\|^p)^{\frac1p}.
\]
Concavity of an operator from a Banach lattice to a Banach space is
defined in a similar manner.

\section{Some basic results on extensions of operators}
\label{sec1}
\setcounter{equation}{0}

In this section we shall prove some general results on extensions of
operators which will be useful for us in the sequel. We start with the
following localization theorem:

\begin{theorem}
\label{thm1.1}
Let $X$ and $Y$ be Banach spaces. Consider the statements:
\begin{itemize}
\item[(i)] Every bounded operator from a arbitrary subspace of $X$ into $Y$
extends to a bounded operator from $X$ to $Y$.
\item[(ii)] There is a constant $K\ge 1$ so that whenever $E\subseteq X$
is a finite dimensional subspace every $T\in B(E,Y)$ admits an extension
$\widetilde{T}\in B(X,Y)$ with $\|\widetilde{T}\|\le K\|T\|$.
\end{itemize}
Then (i) implies (ii) and if $Y$ is a dual space, (ii) implies (i).
\end{theorem}

\bproof
Assume first that (ii) does not hold. By induction we shall construct a
sequence $(E_n)$ of finite dimensional subspaces of $X$, a sequence
$(F_n)$ of subspaces of $X$ of finite codimension and a sequence
$(T_n)\subseteq B(E_n,Y)$ with $\|T_n\|=1$ for all $n\in\bN$ so that the
following conditions are satisfied:
\begin{itemize}
\item[(a)] $F_n\cap \Span\{E_j\mid 1\le j\le n\}=\{0\}$ and the natural
projection of \mbox{$\Span\{E_j\mid 1\le j\le n\}\oplus F_n$} onto
$\Span\{E_j\mid 1\le j\le n\}$ has norm less than or equal to 2 for all
$n\in\bN$.
\item[(b)] $F_{n+1}\subseteq F_n$ for all $n\in\bN$.
\item[(c)] If $\widetilde{T}_n\in B(X,Y)$ is an extension of $T_n$, then
$\|\widetilde{T}_1\|\ge 4$ and $\|\widetilde{T}_n\|\ge 2^{2n+1} \codim
F_{n-1}+\codim F_{n-1}$ for all $n\ge 2$.
\end{itemize}
Since (ii) does not hold, we can for $n=1$ choose a finite dimensional
subspace $E_1$ of $X$ and a $T_1\in B(E_1,Y)$ with $\|T_1\|=1$ so that
any bounded extension of $T_1$ to $X$ has norm greater than or equal to
4. Let $F_1$ be a subspace of finite codimension so that $F_1^\perp$ is
2-norming over $E_1$ ($F_1$ can be chosen to be of codimension $5^{\dim
E_1}$). Clearly $E_1\cap F_1=\{0\}$ and the natural projection of
$E_1\oplus F_1$ onto $E_1$ has norm less than or equal to 2.

Assume now that $E_1,E_2,\dots,E_n$, $F_1,F_2,\dots,F_n$ and
$T_1,T_2,\dots,T_n$ have been constructed so that (a), (b) and (c)
hold. By assumption there is a finite dimensional subspace
$E_{n+1}\subseteq X$ and an operator $T_{n+1}\in B(E_{n+1},Y)$ with
$\|T_{n+1}\|=1$ so that if $\widetilde{T}_{n+1}\in B(X,Y)$ is an
extension of $T_{n+1}$, then
\begin{equation}
\label{eq1.1}
\|\widetilde{T}_{n+1}\|\ge 2^{2n+2} \codim F_n+\codim F_n
\end{equation}
which shows that (c) holds. If we choose a subspace
$\hat{F}_{n+1}\subseteq X$ so that $\hat{F}_{n+1}^\perp$ is 2-norming
over $\Span\{E_j\mid 1\le j\le n\}$ and put $F_{n+1}=\hat{F}_{n+1}\cap
F_n$, then clearly also (a) and (b) are satisfied.

Hence we have constructed the required sequences. Put now $G_1=E_1$ and
$G_{n+1}=E_{n+1}\cap F_n$ for all $n\ge 1$. By choosing an Auerbach
basis for $E_n/G_n$ we easily achieve that there is a subspace
$H_n\subseteq E_n$ and a projection $P_n$ of $X$ onto $H_n$ so that
\begin{eqnarray}
\label{eq1.2}
E_n &=& G_n\oplus H_n\quad\mbox{for all $n\in\bN$}\\
\label{eq1.3}
P_nx&=& 0\quad\mbox{for all $x\in G_n$ and all $n\in\bN$}\\
\label{eq1.4}
\|P_{n+1}\|&\le & \codim F_n\quad\mbox{for all $n\in\bN$}.
\end{eqnarray}
Let $n\ge 2$ and assume that $\widetilde{S}_n\in B(X,Y)$ is an extension
of $T_{n|G_n}$. Put
\begin{equation}
\label{eq1.5}
\widetilde{T}_n = \widetilde{S}_n(I-P_n)+T_nP_n.
\end{equation}
If $x\in E_n$, then
\begin{equation}
\label{eq1.6}
\widetilde{T}_n x = \widetilde{S}_n(x-P_nx)+T_nP_nx = T_n(x-P_nx)+T_nP_n=T_nx.
\end{equation}
Hence $\widetilde{T}_n$ is an extension of $T_n$ and therefore by (c)
\begin{equation}
\label{eq1.6b}
\|\widetilde{T}_n\| \ge 2^{2n+1} \codim F_{n-1}+\codim F_{n-1}
\end{equation}
which in view of \eqref{eq1.4} clearly implies that
\begin{equation}
\label{eq1.7}
\|\widetilde{S}_n\| \ge 2^{2n}.
\end{equation}
By construction $(G_n)$ forms an infinite direct sum and we can
therefore put
\begin{equation}
\label{eq1.8}
G = \bigoplus^\infty_{n=1} G_n.
\end{equation}
We define $S\in B(G,Y)$ by
\begin{equation}
\label{eq1.9}
Sx = \sum^\infty_{n=1} 2^{-n} T_nx_n
\end{equation}
for all $x\in G$ with
\begin{equation}
\label{eq1.10}
x = \sum^\infty_{n=1} x_n \quad x_n\in G_n\quad\mbox{for all $n\in\bN$}.
\end{equation}
(Actually $\|S\|\le 3$).

$S$ does not have a bounded extension to $X$. Indeed, if
$\widetilde{S}\in B(X,Y)$ is an extension, then $2^n\widetilde{S}$ is an
extension of $T_{n|G_n}$ and therefore by \eqref{eq1.7}
\begin{equation}
\label{eq1.11}
\|\widetilde{S}\|\ge 2^n\quad\mbox{for all $n\ge 2$}
\end{equation}
which is a contradiction. This shows that (i) implies (ii).

Assume next that (ii) holds and that $Y$ is a dual space; let $Z$ be a
Banach space so that $Z^*=Y$. Further, let $F\subseteq X$ be a subspace
and $T\in B(F,Z^*)$ with $\|T\|=1$. For every finite dimensional
subspace $E\subseteq F$ we can by assumption find
$\widetilde{T}_E=B(X,Z^*)$ so that
\begin{equation}
\label{eq1.12}
\widetilde{T}_Ex = Tx\quad\mbox{for all $x\in E$,
$\|\widetilde{T}_E\|\le K$}.
\end{equation}
By $\w^*$-compactness it follows that we can find a subnet
$(\widetilde{T}_{E'})$ of $(\widetilde{T}_E)$ and an operator
\mbox{$\widetilde{T}\in B(X,Z^*)$} so that
\begin{equation}
\label{eq1.13}
\widetilde{T}_{E'}x\stackrel{\w^*}{\longrightarrow}
\widetilde{T}x\quad\mbox{for all $x\in X$}.
\end{equation}
Clearly $\widetilde{T}$ is an extension of $T$.
\eproof

The following corollary is an immediate consequence of Theorem
\ref{thm1.1}

\begin{corollary}
\label{cor1.1a}
Let $X$, $Y$ and $Z$ be Banach spaces and assume that $Z$ is finitely
representable in $X$. If every bounded operator from an arbitrary
subspace of $X$ to $Y^*$ extends to a bounded operator from the whole
space to $Y^*$, then every bounded operator from an arbitrary subspace
of $Z$ to $Y^*$ extends.
\end{corollary}
 
Our next result shows that under certain conditions it is enough to
consider extensions of finite rank operators.

\begin{theorem}
\label{thm1.2}
Let $X$ and $Y$ be Banach spaces and $E\subseteq X$ a subspace. Assume
that there is a constant $K$ so that every $T\in E^*\otimes Y$ admits an
extension $\widetilde{T}\in B(X,Y)$ with $\|\widetilde{T}\|\le K\|T\|$.

If either $E$ or $Y$ has the $\lambda$-bounded approximation property,
then every $T\in B(E,Y)$ admits an extension $\widetilde{T}\in
B(X,Y^{**})$ with $\|\widetilde{T}\|\le K\lambda\|T\|$.
\end{theorem}

\bproof
Let $T\in B(E,Y)$. By assumption we can find a net
$(T_\alpha)_{\alpha\in J}\subseteq E^*\otimes Y$ with
$\|T_\alpha\|\le\lambda\|T\|$ for all $\alpha$ so that $T_\alpha x\to
Tx$ for all $x\in E$. Let $\widetilde{T}_\alpha\in B(X,Y)$ denote an 
extension of $T_\alpha$ for each $\alpha\in J$ with
\begin{equation}
\label{eq1.14}
\|\widetilde{T}_\alpha\| \le K\|T_\alpha\|\le K\lambda\|T\|.
\end{equation}
\eqref{eq1.14} immediately gives that there is a $\widetilde{T}\in
B(X,Y^{**})$ with $\|\widetilde{T}\|\le K\lambda\|T\|$ and a subnet
$(\widetilde{T}_{\alpha'})$ of $(\widetilde{T}_\alpha)$ so that
\begin{equation}
\label{eq1.15}
\widetilde{T}_{\alpha'}x\stackrel{\w^*}{\longrightarrow}
\widetilde{T}x\quad\mbox{for all $x\in X$}.
\end{equation}
Since clearly also
$\widetilde{T}_{\alpha'}x\stackrel{\w^*}{\longrightarrow} Tx$ for all
$x\in E$, it follows that $\widetilde{T}$ is the required extension.
\eproof

We shall need:

\begin{lemma}
\label{lemma1.3}
If $E$ is an $n$-dimensional subspace of a Banach space $X$, then
$(E\oplus\ell^n_2)_\infty$ is 12-isomorphic to a subspace of $X$.
\end{lemma}

\bproof
Let $F$ be a subspace of $X$ of finite codimension so that $F^\perp$ is
2-norming on $E$ ($F$ can be chosen so that $\codim F=5^n$). By
Dvoretzky's theorem $F$ contains an $n$-dimensional subspace $G$ with
$d(G,\ell^n _2)\le 2$ and clearly $E\cap G=\{0\}$. It is readily verified that
$(E\oplus\ell^n_2)_\infty$ is 12-isomorphic to $E\oplus G$.
\eproof

The next result shall be very useful for us in the sequel

\begin{theorem}
\label{thm1.4}
Let $X$ and $Y$ be Banach spaces and $\mu$ a measure. If every bounded
operator from an arbitrary subspace of $X$ to $Y^*$ extends to a bounded
operator from $X$ to $Y^*$, then the same holds for every bounded
operator from an arbitrary subspace of $X\oplus L_2(\mu)$ to $Y^*$.
\end{theorem}

\bproof
Let $E\subseteq (X\oplus L_2(\mu))_\infty$ be an arbitrary finite dimensional
subspace. Clearly there exists an $n\in\bN$ so that we can find
$n$-dimensional subspaces $G\subseteq X$ and $F\subseteq L_2(\mu)$ with 
$E\subseteq G\oplus F$. By Lemma \ref{lemma1.3} $G\oplus F$ and
therefore  also $E$ is 12-isomorphic to a subspace of $X$. Hence
$X\oplus L_2(\mu)$ is finitely representable in $X$ and the conclusion
follows from Corollary \ref{cor1.1a}.
\eproof

Finally we shall need the following proposition, the proof of which
is obvious:

\begin{proposition}
\label{prop1.3}
Let $X$ and $Y$ be Banach spaces so that for every subspace  $E\subseteq
X$ every $T\in B(E,Y)$ admits an extension $\widetilde{T}\in B(X,Y)$. If
$Z$ is a quotient of $X$, then $Z$ has the same property.
\end{proposition}

\section{The main results}
\label{sec2}
\setcounter{equation}{0}

We start with the following definition:

\begin{definition}
\label{def2.1}
\begin{itemize}
\item[(i)] A Banach space $X$ is said to have the Maurey extension property ($MEP$) if for any
subspace $E\subseteq X$, any Banach space $Y$ of cotype 2 and every $T\in
B(E,Y)$ there exists an extension $\widetilde{T}\in B(X,Y)$ of $T$.
\item[(ii)] $X$ is said to have $M_p$, $1\le p\le \infty$, if the condition
in (i) holds with $Y=\ell_p$.
\end{itemize}
\end{definition}

Maurey \cite{M1} proved that if $X$ is a Banach space of type 2, then it
has $MEP$. It is readily seen that if a Banach space $X$ has $MEP$, then to
every $\lambda \ge 1$ there exists a constant $C(\lambda)\ge 1$ so that
every bounded operator $T$ from an arbitrary subspace of $X$ to an
arbitrary Banach space $Y$ of cotype $\lambda$ admits an extension
$\widetilde{T}$ from $X$ to $Y$ with $\|\widetilde{T}\|\le C(\lambda)\|T\|$. 

It follows immediately from Theorem \ref{thm1.1} that $X$ has $M_p$ if
and only if there is a constant $K$ so that for every finite dimensional
subspace $E\subseteq X$ every $T\in B(E,\ell_p)$ has an extension
$\widetilde{T}\in B(X,\ell_p)$ with $\|\widetilde{T}\|\le K\|T\|$. We
let $M_p(X)$ denote the smallest constant which can be used here.

Using the above together with the local properties of $L_p$-spaces we
obtain that in Definition \ref{def2.1} we can substitute $\ell_p$ with
an arbitrary infinite dimensional $L_p$-space.

The following result follows immediately from \cite[Theorem 10]{MIP}:

\begin{theorem}
\label{thm3.10}
If $X$ is a Banach space with $M_2$, then it is of weak type 2.
\end{theorem}

We shall postpone the investigation of the property $M_p$ to the next
section and turn to our main results. They state in short that $MEP$
characterizes type 2 spaces among Banach spaces with the Gaussian
average property and that $M_p$, $2<p< \infty$, characterizes Hilbert
spaces among K\"othe function spaces on $[0,1]$. Before we can prove it we need to define certain
special spaces of cotype 2. 

If $\mu$ is a probability measure and $0<\delta<1$, then we define the
space $L_1(\mu;\delta L_2)$ by
\[
L_1(\mu,\delta L_2)=\{(f,\delta f)\mid f\in L_2(\mu)\}\subseteq
(L_1(\mu)\oplus L_2(\mu))_\infty.
\]
Since $L_1(\mu)\oplus L_2(\mu)$ is isomorphic to a subspace of an
$L_1$-space, it follows that $L_1(\mu;\delta L_2)$ is of cotype 2 with a
constant $C$ independent of $\delta$. Note also that it is a sublattice
of $L_1(\mu)\oplus L_2(\mu)$. It is a reflexive space since it is
$\frac{1}{\delta}$-isomorphic to a Hilbert space.

We are now ready to prove:

\begin{theorem}
\label{thm2.1}
If $X$ is a Banach space with the Maurey extension property, then
$B(\ell_{\infty},X^*) = \Pi_2 (\ell_{\infty},X^*)$.
\end{theorem}

\bproof
Let $X$ be a Banach space with $MEP$ and let $(\Omega,\cS,\nu)$ be an
arbitrary probability space. It is clearly enough to show that
$B(X,L_1(\nu)) = \Gamma_2(X,L_1(\nu))$ so let $T\in B(X,L_1(\nu))$ be
arbitrary with $\|T\| = 1$. From \cite[Corollary 1.d.12]{LT2} it follows that if we prove that T is
a 2-convex operator, then we are done. Hence let $n\in\bN$ and
$\{x_1,x_2,\dots,x_n\}\subseteq X$ with $h=\big(\sum^n_{j=1}
|Tx_j|^2\big)^\frac12 \ne 0$. We may assume that $\|h\|_1 =1$. Put
$E=\Span\{x_1,x_2,\dots,x_n\}$, let $\Delta = \{t\in\Omega\mid h(t)>0\}$
and define the probability measure $\mu$ on $\Delta$ by $d\mu =
hd\nu$. Further we let $M_h\colon L_1(\Delta,\nu) \to L_1(\mu)$ denote
the isometry given by:
\begin{equation}
\label{eq2.1}
M_h(f) = fh^{-1} \quad\mbox{for all $f\in L_1(\Delta,\nu)$}
\end{equation} 
and define $\Phi \colon E \to L_1(\mu)$ by $\Phi = M_h T$.
 
Since $X$ has $MEP$ and $L_1(\mu;\delta L_2)$, $0<\delta <1$, has cotype 2
with constant $C$ it follows from Theorem \ref{thm1.4} that there is a
constant $M$ independent of $\delta$ and $\mu$ so that every bounded operator $S$
from a subspace of $(X\oplus L_2(\mu))_{\infty}$ to $L_1(\mu;\delta L_2)$ has an
extension $\widetilde{S}$ to $(X\oplus L_2(\mu))_{\infty}$ with
$\|\widetilde{S}\|\le M\|S\|$. Choose now $\delta$ so that $4CM\delta <1$ and
let $Z\subseteq (X\oplus L_2(\mu))_{\infty}$ be defined by
\begin{equation}
\label{eq2.2}
Z = \{(x,\delta\Phi(x))\mid x\in E\},
\end{equation}
define $I\colon Z\to L_1(\mu;\delta L_2)$ by
\begin{equation}
\label{eq2.3}
I(x,\delta\Phi(x)) = (\Phi(x),\delta\Phi(x))\quad\mbox{for all $x\in E$}
\end{equation}
and let $\widetilde{I}\colon (X\oplus L_2(\mu))_{\infty}\to L_1(\mu;\delta L_2)$ be
an extension of $I$ with $\|\widetilde{I}\|\le M\|I\|\le 2M$. For every
$x\in E$ we now get
\begin{equation}
\label{eq2.4}
(\Phi(x),\delta\Phi(x)) =
\widetilde{I}(x,0)+\delta\widetilde{I}(0,\Phi(x)).
\end{equation}
Using this on the $x_j$'s we obtain
\begin{eqnarray}
\label{eq2.5}
(1,\delta) &=&
\big(\big(\sum^n_{j=1}|\Phi(x_j)|^2\big)^\frac12,\delta\big(\sum^n_{j=1}|\Phi(x_j)|^2\big)^\frac12\big)\\
&=& \big(\sum^n_{j=1}|(\Phi(x_j),\delta\Phi(x_j)|^2\big)^\frac12\\
\nonumber
&=& \big(\sum^n_{j=1} |\widetilde{I}(x_j,0)+\delta
\widetilde{I}(0,\Phi(x_j)|^2\big)^\frac12\\ \nonumber
&\le & \big(\sum^n_{j=1} |\widetilde{I}(x_j,0)|^2\big)^\frac12 +
\delta\big(\sum^n_{j=1}|\widetilde{I}(0,\Phi(x_j)|^2\big)^\frac12.
\end{eqnarray}
Taking norms on both sides of \eqref{eq2.5} we get
\begin{eqnarray}
\label{eq2.6}
1 &\le & \big\|\big(\sum^n_{j=1}
|\widetilde{I}(x_j,0)|^2\big)^\frac12\big\| +
\delta\big\|\big(\sum^n_{j=1}|\widetilde{I}(0,\Phi(x_j)|^2\big)^\frac12\big\|\\
\nonumber
&\le & \big\|\big(\sum^n_{j=1}
|\widetilde{I}(x_j,0)|^2\big)^\frac12\big\| +
\delta C \big(\int^1_0 \big\|\sum^n_{j=1}
r_j(t)\widetilde{I}(0,\Phi(x_j))\big\|^2dt\big)^\frac12\\
\nonumber
&\le &
\big\|\big(\sum^n_{j=1}|\widetilde{I}(x_j,0)|^2\big)^\frac12\big\|+
2\delta CM\big(\int^1_0\big\|\sum^n_{j=1} r_j(t)(0,\Phi(x_j))\big\|^2
dt)\big)^\frac12\\ \nonumber
&=& \big\| \big(\sum^n_{j=1}
|\widetilde{I}(x_j,0)|^2\big)^\frac12\big\| + 2\delta
CM\big\|\big(0,\sum^n_{j=1}|\Phi(x_j)|^2\big)^\frac12\big\| \\ \nonumber
&=& \big\|\big(\sum^n_{j=1} |\widetilde{I}(x_j,0)|^2\big)^\frac12\big\| +
2\delta CM.
\end{eqnarray}
Hence
\begin{equation}
\label{eq2.7}
\frac12 \le \big\|\big(\sum^n_{j=1}
|\widetilde{I}(x_j,0)|^2\big)^\frac12\big\|.
\end{equation}
Let now $Q\colon L_1(\mu)\oplus L_2(\mu)\to L_2(\mu)$ be the canonical
projection onto the second coordinate. By the definition of the order in
$L_1(\mu)\oplus L_2(\mu)$ we have
\[
\big(\sum^n_{j=1} |Q\widetilde{I}(x_j,0)|^2\big)^\frac12 =
Q\big(\sum^n_{j=1} |\widetilde{I}(x_j,0)|^2\big)^\frac12.
\]
Assume now that
\begin{equation}
\label{eq2.8}
\big(\sum^n_{j=1} |\widetilde{I}(x_j,0)|^2\big)^\frac12 = (g,\delta g)
\quad\mbox{with $g\in L_2(\mu)$.}
\end{equation}
If $\big\|\big(\sum^n_{j=1}
|\widetilde{I}(x_j,0)|^2\big)^\frac12\big\|=\|g\|_1$, then by
\eqref{eq2.7}
\begin{eqnarray}
\label{eq2.9}
\frac{\delta}{2} &\le & \delta\big\|\big(\sum^n_{j=1}
|\widetilde{I}(x_j,0)|^2\big)^\frac12\big\| = \delta\|g\|_1\\ \nonumber
&\le & \delta\|g\|_2 = \big\|\big(\sum^n_{j=1} |Q\widetilde{I}(x_j,0)|^2\big)^\frac12\big\|
\end{eqnarray}
and if
$\big\|\big(\sum^n_{j=1}|\widetilde{I}(x_j,0)|^2\big)^\frac12\big\|=\delta\|g\|_2$,
then
\begin{equation}
\label{eq2.10}
\frac12 \le \big\|\big(\sum^n_{j=1}
|\widetilde{I}(x_j,0)|^2\big)^\frac12\big\| = \big\|\big(\sum^n_{j=1} |Q\widetilde{I}(x_j,0)^2\big)^\frac12\big\|.
\end{equation}
Using that the range of $Q\widetilde{I}$ is a Hilbert space we obtain
\begin{equation}
\label{eq2.10a}
\frac{\delta}{2} \le \big\|\big(\sum^n_{j=1}
|Q\widetilde{I}(x_j,0)|^2\big)^\frac12 \big\| = \big(\sum^n_{j=1}
\|Q\widetilde{I}(x_j,0)\|^2\big)^\frac12 \le 2M\big(\sum^n_{j=1}
\|x_j\|^2\big)^\frac12.
\end{equation}
We have now verified that $T$ is 2-convex with constant less than or equal to $4M\delta^{-1}$.
\eproof

Theorem \ref{thm2.1} immediately implies:

\begin{theorem}
\label{thm2.2}
Let $X$ be a Banach space which satisfies one of the following
conditions:
\begin{itemize}
\item[(i)] $X$ has the Gaussian average property.
\item[(ii)] $X$ has the Gordon-Lewis property $\GL_2$ (in particular $X$ could be a Banach lattice).
\item[(iii)] $X$ is isomorphic to a subspace of a Banach lattice of
finite cotype.
\end{itemize}
If $X$ has the Maurey extension property, then $X$ is of type 2.
\end{theorem}

\bproof
Let $X$ be a Banach space with $MEP$.
\begin{itemize}
\item[(i)] If $X$ has $\GAP$, then it follows from Theorem \ref{thm2.1}
and \cite[Theorem 1.10]{CN1} that $X$ is of type 2.
\item[(ii)] Since $X$ has $MEP$, it is of finite cotype and if in
addition it has $\GL_2$, then it has $\GAP$ by \cite[Theorem 1.3]{CN1}. (ii)
can also be derived directly from Theorem \ref{thm2.1} and
\cite[Proposition 8.16]{P1}.
\item[(iii)] If $X$ is isomorphic to a subspace of a Banach lattice of
finite cotype, then it has $\GAP$ by \cite[Theorem 1.4]{CN1}.
\end{itemize}
\eproof

\bremark
It follows from \cite{CN1} that every space of type 2 has $\GAP$. Hence if
there exists a Banach space with $MEP$ and without $\GAP$, then it cannot
have type 2.

If a Banach space $X$ has $MEP$, then every bounded operator from a
subspace of $X$ to a cotype 2 space $Y$ with $\GL$ can be extended to
$X$ through a Hilbert space (as in Maurey's original result). Indeed, let
$E$ be a subspace of $X$ and $T\in B(X,Y)$. Since $E$ has $MEP$ and $Y$ has
$\GL(1,2)$ by \cite[Theorem 3.4]{CN2}, it follows from Theorem
\ref{thm2.1} and Theorem \ref{thm3.1} in the next section that $T\in
\Gamma_2(E,Y)$. Since $X$ has $MEP$, the part of the factorization of $T$
which goes into a Hilbert space can be extended to $X$.
\eremark

Before we can prove our main result on $M_p$, $2<p<\infty$, we need a
sequence space equivalent of the spaces considered in Theorem \ref{thm2.1}. 

If $X$, respectively $Y$, have unconditional normalized bases $(x_n)$,
respectively $(y_n)$, then we say that $(x_n)$ dominates $(y_n)$ and
write $(y_n)<(x_n)$ if the linear operator
$T\colon\Span(x_n)\to\Span(y_n)$ defined by $Tx_n=y_n$ for all $n\in\bN$
is bounded. If $1\le q\le\infty$ and the unit vector basis of $\ell_q$
dominates $(x_n)$, respectively is dominated by $(x_n)$, then we shall
say that $(x_n)$ satisfies an upper $p$-estimate, respectively
lower $p$-estimate.

If $1\le q<\infty$ and $(e_n)$ denotes the unit vector basis of
$\ell_q$, then for every $0<\delta<1$ we define the space $X(\delta,q)$
to be the closed linear span in $(X\oplus \ell_q)_\infty$ of the
sequence $(x_j+\delta e_j)$.

The next theorem which shall be very useful for us in several contexts states:
\begin{theorem}
\label{thm2.3}
Let $X$, respectively $Y$, be Banach spaces with normalized
unconditional bases $(x_n)$, respectively $(y_n)$, $1\le q<\infty$, so
that $(y_n)<(x_n)$ with constant $K_1$ and $(y_n)$ satisfies an upper
$q$-estimate with constant $K_2$. If for some $0<\delta<1$ the formal
identity operator $I_\delta$ from $X(\delta,q)$ to $Y(\delta,q)$ extends
to a bounded operator $\widetilde{I}_\delta$ from
$(X\oplus\ell_q)_\infty$ to $Y(\delta,q)$ with
$\|\widetilde{I}_\delta\|<\delta^{-1}$, then for all $(t_n)\subseteq\bR$
\begin{eqnarray}
\label{eq2.11}
\delta^2(1-\|I_\delta\|\delta)\big(\sum^\infty_{n=1}|t_n|^2\big)^{\frac12}
&\le & \sqrt{2} K_2 ubc(x_n)\big\|\sum^\infty_{n=1}
t_nx_n\big\|\quad\mbox{if $1\le q\le 2$}\\
\label{eq2.12}
\delta^2(1-\|I_\delta\|\delta)\big(\sum^\infty_{n=1}|t_n|^q\big)^{\frac1 q}
&\le & K_2 ubc(x_n)\big\|\sum^\infty_{n=1}
t_nx_n\big\|\quad\mbox{if $2\le q\le \infty$}
\end{eqnarray}
e.g. $(x_n)$ has a lower 2-estimate if $1\le q\le 2$ and a lower
$p$-estimate if $2\le q<\infty$.
\end{theorem}

\bproof
Since $\widetilde{I}_\delta$ extends $I_\delta$, we have for all
$n\in\bN$
\begin{equation}
\label{eq2.13}
y_n+\delta e_n = \widetilde{I}_\delta x_n+\delta\widetilde{I}_\delta e_n
\end{equation}
and hence by the triangle inequality
\begin{equation}
\label{eq2.14}
(1-\|\widetilde{I}_\delta\|\delta) \le \|\widetilde{I}_\delta
x_n\|\quad\mbox{for all $n\in\bN$}.
\end{equation}
Let $Q\colon (Y\oplus \ell_q)_\infty\to \ell_q$ be the canonical
projection and let $T=Q\widetilde{I}_\delta$. Fix $n\in\bN$ and let
$(a_k)\subseteq\bR$ be chosen so that
\begin{equation}
\label{eq2.15}
\widetilde{I}_\delta x_n=\sum^\infty_{k=1} a_ky_k + \delta
\sum^\infty_{k=1} a_ke_k.
\end{equation}
If $\|\widetilde{I}_\delta
x_n\|=\delta\big(\sum^\infty_{k=1}|a_k|^q\big)^{\frac1q}$, then by
\eqref{eq2.14}
\begin{equation}
\label{eq2.16}
(1-\|\widetilde{I}_\delta\|\delta) \le
\delta\big(\sum^\infty_{k=1}|a_k|^q\big)^{\frac1q}= \|T x_n\|
\end{equation}
and if $\|\widetilde{I}_\delta x_n\| = \big\|\sum^\infty_{k=1}
a_ky_k\big\|$, we obtain
\begin{equation}
\label{eq2.17}
\delta(1-\|\widetilde{I}_\delta\|\delta) \le
\delta\big\|\sum^\infty_{k=1} a_ky_k\big\|\le
K_2\delta\big(\sum^\infty_{k=1} |a_k|^q\big)^{\frac1q}=\|T x_n\|.
\end{equation}
Comparing \eqref{eq2.16} and \eqref{eq2.17} we get that for all $n\in\bN$
\begin{equation}
\label{eq2.18}
K_2^{-1}\delta(1-\|\widetilde{I}_\delta\|\delta)\le\|T x_n\|.
\end{equation}
Let $r=\max(q,2)$. Since $\ell_q$ is of cotype $r$, we get for all
$n\in\bN$ and all $(t_j)^n_{j=1}\subseteq\bR$:
\begin{eqnarray}
\label{eq2.19}
K_2^{-1}\delta(1-\|\widetilde{I}\|\delta)\big(\sum^n_{j=1}
|t_j|^r\big)^{\frac1r} &\le & \big(\sum^n_{j=1} |t_j|^r \|T
x_j\|^r\big)^{\frac1r}\\ \nonumber
&\le & C_q\big(\int^1_0 \big\|\sum^n_{j=1} r_j(t)t_j Tx_j\big\|^r
dt\big)^{\frac1r} \\ \nonumber
&\le & C_q\|T\|\big(\int^1_0 \big\|\sum^n_{j=1} r_j(t)t_jx_j\big\|^r
dt\big)^{\frac1r}\\ \nonumber
&\le & C_q\delta^{-1} ubc(x_j)\big\|\sum^n_{j=1} t_jx_j\big\|
\end{eqnarray}
where $C_q\le \sqrt{2}$ for $1\le q<2$ and $C_q=2$ for $2\le
p<\infty$. \eqref{eq2.19} immediately gives  \eqref{eq2.11} and
\eqref{eq2.12}. Note that our assumptions imply that $\delta<K_1^{-1}$.
\eproof

\bremark
Theorem \ref{thm2.3} remains true if we assume that both $X$ and $Y$ are
finite dimensional.
\eremark

Theorem \ref{thm2.3} was inspired by Nigel Kalton, who drew our
attention to the spaces $\ell_p(\delta,2)$, $p>2$ in order to prove that
$\ell_p$ does not have $M_r$ for $2<p<r<\infty$ which subsequently lead
to the idea of the proof of Theorem \ref{thm2.1}. Spaces like
$\ell_p(\delta,2)$ were first considered by Rosenthal in his construction
of new ${\cal L}_p$ spaces \cite{R}.

Before we go on we need a few facts about the spaces $\ell_p(\delta,2)$,
$p>2$, which all go back to \cite{R}. Hence let $2<p<\infty$ and
$0<\delta<1$. The space $L_p(0,\infty)\cap L_1(0,\infty)$ equipped with
the maximum of the $p$-norm and the 2-norm is a rearrangement invariant
function space on $[0,\infty[$ which is isomorphic to $L_p(0,1)$,
\cite[Theorem 2.f.1]{LT2}. In addition $\ell_p(\delta,2)$ is isometric
to a norm 1 complemented subspace of $L_p(0,\infty)\cap
L_2(0,\infty)$. Indeed, it is readily seen that if we take a sequence
$(I_k)^\infty_{k=1}$ of mutually disjoint intervals in $[0,\infty[$ each
of length $\delta^{\frac{2p}{n-1}}$, then the closed linear span of
$\{1_{I_k}\}$ is isometric to $\ell_p(\delta,2)$. This span is also norm
1 complemented since conditional expectations are norm 1 projections in
$L_p(0,\infty)\cap L_2(0,\infty)$. Hence we have verified:

\begin{lemma}
\label{lemma2.2}
Let $2<p<\infty$. There exists a constant $C$ so that for all
$\delta\in]0,1[$ $\ell_p(\delta,2)$ is $C$-isomorphic to a
$C$-complemented subspace of $L_p(0,1)$.
\end{lemma}

We need yet another lemma:

\begin{lemma}
\label{lemma2.3}
If $X$ is a Banach space with $M_p$ for some $2<p<\infty$, then
$\inf\{q\mid X$ has cotype $q\}<p$. In particular $X$ has cotype $p$.
\end{lemma}

\bproof
Put $q_0=\inf\{q\mid X$ has cotype $q\}$. By \cite{MP} $L_{q_0}(0,1)$ is
finitely representable in $X$ and hence it has $M_p$ by Corollary
\ref{cor1.1a}. If $p\le q_0$, then $L_p(0,1)$ is a quotient of
$L_{q_0}(0,1)$ and hence it also has $M_p$ by Proposition \ref{prop1.3}; this
is a contradiction since $L_p(0,1)$ contains uncomplemented subspaces
isomorphic to $\ell_p$ \cite{R}.
\eproof

We are now ready to prove:

\begin{theorem}
\label{thm2.4}
If $2<p<\infty$ and $X$ is a Banach space with $M_p$, then the following
statements hold:
\begin{itemize}
\item[(i)] For every $\lambda\ge 1$ there exists a constant $c(\lambda)$
so that whenever $(x_j)\subseteq X$ is a finite or infinite
$\lambda$-unconditional normalized sequence then
\begin{equation}
\label{eq2.20}
c(\lambda)\big(\sum_j |a_j|^2\big)^{\frac12}\le\big\|\sum_j
a_jx_j\big\|\qquad\mbox{for all $(a_j)\subseteq\bR$}.
\end{equation}
\item[(ii)] $X$ is of weak type 2 and has property $(H)$. If in addition
$X$ is a Banach lattice then it is a weak Hilbert space which satisfies
a lower 2-estimate.
\end{itemize}
\end{theorem}

\bproof
\begin{itemize}
\item[(i)] Let $n\in\bN$, $\lambda\ge 1$ and let $(x_j)^n_{j=1}\subseteq
X$ be a normalized $\lambda$-unconditional sequence. Since $([x_j]\oplus
\ell^n_2)_\infty$ is 12-isomorphic to a subspace of $X$, it follows that
$([x_j]\oplus\ell_2^n)_\infty$ has $M_p$ with constant less than or
equal to $12M_p(X)$. Combining this with Lemma \ref{lemma2.2} we get
that every bounded operator $T$ from a subspace of
$([x_j]\oplus\ell^n_2)$ to any $\ell_p(\delta,2)$, $0<\delta<1$, has an
extension $\widetilde{T}$ to $([x_j]\oplus\ell_2^n)_\infty$ with
$\|\widetilde{T}\|\le 12C^2 M_p(X)$. By Lemma \ref{lemma2.3} $X$ has
cotype $p$ and hence the cotype constant of
$([x_j]\oplus\ell^n_2)_\infty$ is less than or equal to $2 K_p(X)$ and
therefore the formal identity operator $I_\delta$ of $[x_j](\delta,2)$
into $\ell_p(\delta,2)$ has a norm less than or equal to $2K_p(X)$. If we now
choose $\delta$ so that $24C^2k_p(X)M_p(X)\delta <1$, then it follows
that $I_\delta$ has an extension to $([x_j]\oplus\ell^n_2)_\infty$ with
norm less than $\delta^{-1}$. Hence by Theorem \ref{thm2.1}
we get for all $(t_j)^n_{j=1}\subseteq\bR$:
\begin{equation}
\label{eq2.21}
\frac{\delta^2}{2} \big(\sum^n_{j=1}
|t_j|^2\big)^{\frac12}\le\lambda\big\|\sum^n_{j=1} t_jx_j\big\|
\end{equation}
which proves \eqref{eq2.20}.
\item[(ii)] Since $X$ has $M_p$, it also has $M_2$ (because $L_p$ has a
complemented subspace isomorphic to a Hilbert space) and hence $X$ is of
weak type 2. Combining this with \eqref{eq2.20} we get that there exists
a constant $C(\lambda)$ so that if $(x_j)^n_{j=1}\subseteq X$ is
$\lambda$-unconditional and normalized, then
\begin{equation}
\label{eq2.22}
c(\lambda)\sqrt{n}\le\big\|\sum^n_{j=1} x_j\big\|\le C(\lambda)\sqrt{n},
\end{equation}
which proves that $X$ has property $(H)$.

If in addition $X$ is a Banach lattice, then it follows from
\cite[Corollary 4.4]{NTJ} that $X$ is a weak Hilbert space which by
\eqref{eq2.20} satisfies a lower 2-estimate.
\end{itemize}
\eproof

Let us conclude this section with two corollaries.

\begin{corollary}
\label{cor2.5}
Let $X$ be a Köthe function space on $[0,1]$. If $X$ has $M_p$ for some
$p$, $2<p<\infty$, then $X$ is lattice isomorphic to $L_2(0,1)$.
\end{corollary}

\bproof
It follows from theorem \ref{thm2.4} that $X$ is a weak Hilbert space
and hence by \cite[Theorem 3]{NJN2} $X$ is lattice isomorphic to $L_2(0,1)$.
\eproof

\begin{corollary}
\label{cor2.6}
If $X$ is a Banach lattice with an upper 2-estimate which has $M_p$ for
some $p$, $2<p<\infty$, then $X$ is isomorphic to a Hilbert space.
\end{corollary}

\section{The extension properties $M_p$, $1\le p<\infty$}
\label{sec3}

In this section we shall investigate the properties $M_p$ in greater
detail. Our first theorem gives a necessary and sufficient condition for an
operator from a subspace of $X$ to $\ell_p$ to be extended to $X$.

\begin{theorem}
\label{thm3.4}
Let $X$ be a Banach space, $E$ a subspace of $X$ and $T\in B(E,\ell_p)$,
$1\le p\le\infty$. Let $Q$ be the natural quotient map of $X^*$ onto
$E^*$. The following statements are equivalent:
\begin{itemize}
\item[(i)] $T$ has an extension $\widetilde{T}\in B(X,\ell_p)$.
\item[(ii)] There is a constant $K\ge 1$ so that for all Banach spaces
$Z$ and all $S\in B(Z,E)$ with $S^*Q\in\pi_p(X^*,Z^*)$ $TS$ is
$p$-integral with
\begin{equation}
\label{eq3.3}
i_p(TS) \le K\pi_p(S^*Q).
\end{equation}
\end{itemize}
\end{theorem}

\bproof
Assume that (i) holds and let $\widetilde{T}\in B(X,\ell_p)$ be an
extension. Since $\|\widetilde{T}\| = \gamma_p(\widetilde{T})$, it
follows from \cite[Theorem 9.11]{DJT} that if $Z$ is an arbitrary
Banach space and $S\in B(Z,E)$ with $S^*Q\in\pi_p(X^*,Z^*)$, then
$\widetilde{T}S=TS$ is $p$-integral with
\[
i_p(TS) = i_p(\widetilde{T}S) \le \|\widetilde{T}\| \pi_p(S^*Q)
\]
which is \eqref{eq3.3} with $K=\|\widetilde{T}\|$.

Assume next that (ii) holds and define
\begin{equation}
\label{eq3.4}
\cN = \{U\in N_1(\ell_p,X)\mid U(\ell_p)\subseteq E\}.
\end{equation}
If we can prove that $T$ acts as a bounded linear functional on $\cN$ via
trace duality, then since $N_1(\ell_p,X)^*=B(X,\ell_p^{**})$ it follows
that $T$ admits an extension $\widetilde{T}\in B(X,\ell_p)$.

Hence let $U\in \cN$ be arbitrary and let $\e>0$. From Kwapien's
characterization of $\Gamma_p^*$ \cite{K1} it follows that there exist a Banach
space $Z$, $A\in\pi_{p'}(\ell_p,Z)$ and $S\in B(Z,E)$ with $S^*Q
\in\pi_p(X^*,Z^*)$, so that $U=SA$ and
\begin{equation}
\label{eq3.5}
\pi_{p'}(A)\pi_p(S^*Q)\le\nu_1(U)+\e.
\end{equation}
Applying now (\ref{eq1.3}) we obtain
\begin{equation}
\label{eq3.6}
|\tr(TU)| \le i_p(TS)\pi_{p'}(A)\le K\pi_p(S^*Q)\pi_{p'}(A)\le
K(\nu_1(U)+\e).
\end{equation}
Since $\e>0$ was arbitrary, \eqref{eq3.6} shows that $T$ admits an
extension $\widetilde{T}$ with $\|\widetilde{T}\|\le K$.
\eproof

In our next result we shall use Theorem \ref{thm3.4} to give a necessary
and sufficient condition for every operator from a given subspace
of $X$ to extend to $X$.

\begin{theorem}
\label{thm3.5}
Let $E$ be a subspace of a Banach space $X$ and $1\le p\le 2$. Further
let $Q$ denote the canonical quotient map of $X^*$ onto $E^*$. The
following statements are equivalent
\begin{itemize}
\item[(i)] Every $T\in B(E,\ell_p)$ extends to a $\widetilde{T}\in
B(X,\ell_p)$.
\item[(ii)] There is a constant $K\ge 1$ so that every $T\in E^*\otimes
\ell_p$ extends to a $\widetilde{T}\in B(E,\ell_p)$ with
$\|\widetilde{T}\|\le K\|T\|$.
\item[(iii)] There exists a constant $K\ge 1$ so that for all Banach
spaces we have that whenever $S\in B(E^*,Z)$ with $SQ\in\pi_p(E^*,Z)$
then $S\in\pi_p(E^*,Z)$ with
\begin{equation}
\label{eq3.7}
\pi_p(S)\le K\pi_p(SQ).
\end{equation}
\end{itemize}
\end{theorem}

\bproof
In view of the open mapping theorem and Theorem \ref{thm1.2} it is
immediate that (i) and (ii) are equivalent. Hence assume that (ii) holds
and let $K$ be a constant from there. Let $Z$ be an arbitrary Banach space and let $S\in B(E^*,Z)$
with $SQ \in\pi_p(E^*,Z)$.Our assumption and \cite{K2}
(see also \cite{NJN1}) imply that
\begin{eqnarray}
\label{eq3.8}
\sup\{\|TS^*\|_m\mid T & \in & B(E^{**},\ell_p),\|T\|\le 1\}\\ \nonumber
&\le & K
\sup\{\|TS^*\|_m\mid T\in B(X^{**},\ell_p),\|T\|\le 1\}\\ \nonumber
&=& K\pi_p(SQ).
\end{eqnarray}
Since the left hand side is finite, we can conclude that it is equal to
$\pi_p(S)$. Hence $S\in\pi_p(E^*,Z)$ with $\pi_p(S)\le K\pi_p(SQ)$.

Assume next that (iii) holds and let $T\in B(E,\ell_p)$ be arbitrary. We
shall verify that (ii) of Theorem \ref{thm3.4} holds. Hence let $Z$ be
an arbitrary Banach space and $S\in B(Z,E)$ with
$S^*Q\in\pi_p(X^{**},Z^*)$. From (\ref{eq3.7}) we conclude that
$S^*\in\pi_p(E^*,Z^*)$, and therefore by \cite{K2} $TS$ is order
bounded and hence also $p$-integral with
\begin{equation}
\label{eq3.9}
i_p(TS) \le \|TS\|_m\le \|T\| \pi_p(S^*)\le K\|T\|\pi_p(S^*Q).
\end{equation}
Hence $T$ admits an extension $\widetilde{T}$ to $X$ with
$\|\widetilde{T}\|\le K\|T\|$.
\eproof

Using the previous results we now obtain:

\begin{theorem}
\label{thm3.6}
Let $X$ be a Banach space and $1\le p\le\infty$. The following
statements are equivalent.
\begin{itemize}
\item[(i)] $X$ has $M_p$.
\item[(ii)] There exists a constant $K\ge 1$ so that if $E$ is an
arbitrary subspace of $X$, $Q_E$ is the canonical quotient map of  $X^*$ onto
$E^*$ and $Z$ is an arbitrary Banach space, then for every $S\in
B(E^*,Z)$ with $SQ\in\pi_p(X^*,Z)$ we have that $S\in\pi_p(E^*,Z)$ with
\begin{equation}
\label{eq3.10}
\pi_p(S)\le K\pi_p(SQ).
\end{equation}
\end{itemize}
\end{theorem}

\bproof
The equivalence follows immediately from Theorem \ref{thm1.1} and
Theorem \ref{thm3.5}.
\eproof

We now need the following lemma:

\begin{lemma}
\label{lemma3.7}
If $X$ is a Banach space with $M_1$, then there is a $p$, $1<p\le 2$ so that $X$
has type $p$.
\end{lemma}

\bproof

Let $X$ have $M_1$. If $X$ is not of type greater than one, then by \cite{MP}
$\ell_1$ is finitely representable in $X$ and hence it follows from
Corollary \ref{cor1.1a} that $\ell_1$ has
$M_1$. By \cite{B} $\ell_1$ contains an uncomplemented subspace $E$ 
isomorphic to $\ell_1$; hence no isomorphism of $E$ onto $\ell_1$ can be
extended to $\ell_1$ which is a contradiction.
\eproof

We are now able to prove

\begin{theorem}
\label{thm3.8}
If $X$ is a Banach space, then the following statements hold
\begin{itemize}
\item[(i)] If $X$ has $M_1$, then it has $M_2$.
\item[(ii)] If $1<p<2$, then $X$ has $M_1$ if and only if it has $M_p$.
\item[(iii)] If $X$ has $M_p$ for some $p$, $2<p<\infty$ then it has $M_2$.
\end{itemize}
\end{theorem}

\bproof
\begin{itemize}
\item[(i)] Let $X$ have $M_1$. By Lemma \ref{lemma3.7} there is a $q>1$
so that $X$ has type $q$ and let $1<p<q$. If $E\subseteq X$ is a
subspace, then it follows from \cite{MP} that $\pi_1(E^*,Z)=\pi_p(E^*,Z)$
for every Banach space $Z$ and hence we get from our assumption and
Theorem \ref{thm3.6} that $X$ has $M_p$. Since $L_p(0,1)$ has a
complemented subspace isomorphic to a Hilbert space, we obtain that $X$
has $M_2$.
\item[(ii)] Let $1<p<2$ and assume first that $X$ has $M_1$. By (i) and
Theorem \ref{thm3.10} $X$ has type $q$ for all $q<2$ and hence we can
argue like in (i) to get that $X$ has $M_p$. Assume next that $X$ has
$M_p$. Again the argument of (i) shows that $X$ has $M_2$ and is
therefore of type $q$ for all $q<2$. If $E\subseteq X$ is a subspace and
$T\in B(E,\ell_1)$, then $T\in\Gamma_p(E,\ell_1)$ and hence it can be
extended to a bounded  $\widetilde{T}\in B(X,\ell_1)$.
\item[(iii)] If $2<p<\infty$, then $L_p(0,1)$ has a complemented subspace
isomorphic to a Hilbert space and hence if $X$ has $M_p$, it also has
$M_2$.
\end{itemize}
\eproof

We shall now need the following factorization theorem which is a
generalization of \cite[Theorem 8.17]{P1}.

\begin{theorem}
\label{thm3.1}
Let $1\le p \le 2$ and let $X$ and $Y$ be Banach spaces . If
$B(\ell_{\infty},X^*) = \Pi_{p'}(\ell_{\infty},X^*)$ and  $Y$ has
$\GL(1,p)$, then $B(X,Y)\subseteq\Gamma_p(X,Y^{**})$ and
\begin{equation}
\label{eq3.1}
\gamma_p(T) \le C_q(X^*)\GL_{1,p}(Y)\|T\|\qquad\mbox{for all $T\in
B(X,Y)$}.
\end{equation}
\end{theorem}

\bproof
Let $T\in B(X,Y)$ be arbitrary. We shall use \cite[Theorem 9.11]{DJT} to
show that $T\in\Gamma_p(X,Y^{**})$. To this end let $Z$ be an arbitrary Banach
space and $S\in B(Z,X)$ with $S^*\in\pi_p(X^*,Z^*)$. The assumptions on
$X$ give that $S^*$ is absolutely summing and since $Y$ has $\GL(1,p)$, we
get that $TS$ is $p$-integral with
\begin{equation}
\label{eq3.2}
i_p(TS) \le \GL_{1,p}(Y)\pi_1(S^*T^*)\le
C_q(X^*)\GL_{1,p}(Y)\pi_p(S^*)\|T\|.
\end{equation}
\eqref{eq3.2} together with the above-mentioned theorem gives \eqref{eq3.1}.
\eproof

\begin{corollary}
\label{cor3.3}
Let $p,q$ and $X$ be as in Theorem \ref{thm2.3}. If $Y$ is a
complemented subspace of a $p$-concave Banach lattice $Z$, then
$B(X,Y)=\Gamma_p(X,Y)$.
\end{corollary}

\bproof
It follows from \cite{GJN} that $Y$ has $\GL(1,p)$ and since $Z$ does
not contain $c_0$, it follows from \cite{LT2} that $Z$ and hence also $Y$
is complemented in its second dual.
\eproof

The next theorem is a direct consequence of Theorems \ref{thm3.1} and
\ref{thm3.8}.

\begin{theorem}
\label{thm3.9}
Let $X$ be a Banach space with $M_1$ and $Y$ a Banach space with
$\GL(1,p)$ where $1\le p<2$. If $E\subseteq X$ is a subspace, then every
$T\in B(E,Y)$ extends to a $\widetilde{T}\in B(X,Y^{**})$ with
\begin{equation}
\label{eq3.11}
\|\widetilde{T}\| \le M_p(X)\GL_{1,p}(Y) T_r(X)\|T\|\qquad\mbox{for all
$r$, $p<r<2$.} 
\end{equation}
\end{theorem}

\bproof
Choose $p<r<2$ and let $T\in B(E,Y)$. Since $X$ (and hence $E$) has type
$r$ by Theorem \ref{thm3.8}, we get from Theorem \ref{thm3.1} that
$T\in\Gamma_p(E,Y^{**})$ with
\begin{equation}
\label{eq3.12}
\gamma_p(T)\le T_r(X)\GL_{1,p}(Y) \|T\|.
\end{equation}
Since $X$ also has $M_p$ it follows from (\ref{eq3.12}) that $T$ can be
extended to a $\widetilde{T}\in B(X,Y^{**})$ so that (\ref{eq3.11})
holds.
\eproof

It is immediate from the definition of $M_2$ that the following holds:

\begin{proposition}
\label{prop3.11}
Let $X$ be a Banach space with $M_2$. For every finite dimensional
subspace $E\subseteq X$ there exists a projection $P$ of $X$ onto $E$
with
\begin{equation}
\label{eq3.13}
\|P\| \le M_2(X)d(E,\ell_2^{dimE}).
\end{equation}
\end{proposition}

If $X$ is a Banach space and there exists a constant $K$ so that
(\ref{eq3.13}) holds with $K$ interchanged with $M_2(X)$, then $X$ is
said to have the {\em Maurey projection property}. It follows from
\cite[Theorem 11.6]{P1} that a Banach space with this property is of
weak type 2. We end this section with the following result:

\begin{theorem}
\label{thm3.12}
Let $X$ be a K\"othe function space on $[0,1]$ with an unconditional
basis. If $X$ has the Maurey projection property, then it is of type 2.
\end{theorem}

\bproof
Since $X$ has an unconditional basis, it follows from \cite{KW} that $X$
is isomorphic to $X(\ell_2)$ ($=\ell_2 \otimes_m X$). It therefore
follows from from \cite[Remark 11.8]{P2} that $X$ being of weak type 2
is actually of type 2.
\eproof

\vspace{1cm}

\noindent Department of Mathematics,\\ University of Missouri,\\ Columbia MO 65211,\\
pete@casazza.math.missouri.edu\\

\noindent Department of Mathematics and Computer Science,\\ SDU-Odense 
University,\\ Campusvej 55, DK-5230 Odense M, Denmark,\\
njn@imada.sdu.dk

\end{document}